# A simple approach to q- Chebyshev polynomials

Johann Cigler

Fakultät für Mathematik, Universität Wien

*johann.cigler@univie.ac.at*

**Abstract**

It is shown that some $q-$analogues of the Fibonacci and Lucas polynomials lead to $q-$analogues of the Chebyshev polynomials which retain most of their elementary properties.

## 1. Introduction

Let me first sketch the context from which this paper originated. I have been interested in $q-$analogues of bivariate Fibonacci and Lucas polynomials with simple formulae as well as simple recurrences. The Fibonacci polynomials admit many such analogues but for the Lucas polynomials the situation is more complicated (see e.g. [5]). It turned out that by introducing an additional parameter both aims can be accomplished. A special choice of this parameter led me to a $q-$analogue of the bivariate Chebyshev polynomials which admits simple generalizations of many properties of the classical univariate Chebyshev polynomials. Some variants of these polynomials (i.e. the Al-Salam and Ismail polynomials) previously occurred in [1] and have been studied in [6] from another point of view.

To make the paper intelligible for non-specialists let me first recall some well-known facts about these polynomials.

The *Fibonacci polynomials* $F_n(x,s)$ are defined by the recursion

$$F_n(x,s) = xF_{n-1}(x,s) + sF_{n-2}(x,s) \tag{1.1}$$

with initial values $F_0(x,s) = 0$ and $F_1(x,s) = 1.$ They have the explicit expression

$$F_n(x,s) = \sum_{k=0}^{\left\lfloor \frac{n-1}{2} \right\rfloor} \binom{n-1-k}{k} s^k x^{n-1-2k}. \tag{1.2}$$

The *Lucas polynomials* $L_n(x,s)$ satisfy the same recursion

$$L_n(x,s) = xL_{n-1}(x,s) + sL_{n-2}(x,s) \tag{1.3}$$

with initial values $L_0(x,s)=2$ and $L_1(x,s)=x.$ They are given by

$$L_n(x,s)=\sum_{k=0}^{\left\lfloor \frac{n}{2}\right\rfloor}\frac{n}{n-k}\binom{n-k}{k}s^k x^{n-2k}. \tag{1.4}$$

These polynomials satisfy a multitude of identities most of which can be proved by using *Binet's formulae*

$$F_n(x,s)=\frac{\alpha^n-\beta^n}{\alpha-\beta} \tag{1.5}$$

and

$$L_n(x,s)=\alpha^n+\beta^n \tag{1.6}$$

where $\alpha=\frac{x+\sqrt{x^2+4s}}{2}$ and $\beta=\frac{x-\sqrt{x^2+4s}}{2}$ are the roots of the equation $z^2=xz+s.$

Unfortunately there are no simple $q-$analogues of Binet's formulae. Useful substitutes are the *Fibonacci matrices*

$$C^n=\begin{pmatrix} F_{n-1}(x,s) & F_n(x,s) \\ F_n(x,s) & F_{n+1}(x,s)\end{pmatrix} \tag{1.7}$$

where $C=\begin{pmatrix} 0 & 1 \\ s & x\end{pmatrix}$ satisfies the equation $C^2=xC+sI.$

Here $I=\begin{pmatrix} 1 & 0 \\ 0 & 1\end{pmatrix}$ is the $2x2-$identity matrix.

The Lucas polynomial $L_n(x,s)$ coincides with the trace

$$L_n(x,s)=tr\left(C^n\right)=F_{n+1}(x,s)+F_{n-1}(x,s) \tag{1.8}$$

of the Fibonacci matrix $C^n$. The reader is referred to [5] for details.

As is well known a sequence $\left(p_n\right)$ of monic polynomials of degree $\deg p_n=n$ is *orthogonal with respect to some linear functional* $\Lambda,$ i.e. $\Lambda\left(p_n p_m\right)=0$ for $n\neq m$ and $\Lambda\left(p_n^2\right)\neq 0$, if and only if it satisfies a three-term recurrence of the form $p_n(x)=(x-s_n)p_{n-1}(x)-t_n p_{n-2}(x).$ This functional is uniquely determined by $\Lambda\left(p_n\right)=[n=0].$ The numbers $\Lambda\left(x^n\right)$ are called the associated *moments*. The Fibonacci and Lucas polynomials satisfy the 3-term recurrence (1.1) and are therefore orthogonal.

The *bivariate Chebyshev polynomials* $T_n(x,s)$ *of the first kind* can be defined by the recursion

$$T_n(x,s)=2xT_{n-1}(x,s)+sT_{n-2}(x,s) \tag{1.9}$$

with initial values $T_0(x,s)=1$ and $T_1(x,s)=x.$
The *classical Chebyshev polynomials* $T_n(x)=T_n(x,-1)$ are characterized by the identity

$$\cos\left(n\vartheta\right)=T_n(\cos\vartheta). \tag{1.10}$$

The polynomials $T_n(x,s)$ are related to the Lucas polynomials by

$$T_n(x,s)=2^{n-1}L_n\left(x,\frac{s}{4}\right)=\frac{1}{2}\left(\left(x+\sqrt{x^2+s}\right)^n+\left(x-\sqrt{x^2+s}\right)^n\right). \tag{1.11}$$

The *bivariate Chebyshev polynomials* $U_n(x,s)$ *of the second kind* can be defined by the same recursion

$$U_n(x,s)=2xU_{n-1}(x,s)+sU_{n-2}(x,s) \tag{1.12}$$

with initial values $U_0(x,s)=1$ and $U_1(x,s)=2x.$
The classical polynomials $U_n(x)=U_n(x,-1)$ are characterized by the identity

$$\frac{\sin\left((n+1)\vartheta\right)}{\sin\vartheta}=U_n(\cos\vartheta). \tag{1.13}$$

The bivariate polynomials of the second kind are related to the Fibonacci polynomials by

$$U_n(x,s)=2^nF_{n+1}\left(x,\frac{s}{4}\right)=\frac{\left(x+\sqrt{x^2+s}\right)^{n+1}-\left(x-\sqrt{x^2+s}\right)^{n+1}}{2\sqrt{x^2+s}}. \tag{1.14}$$

Both types of Chebyshev polynomials are characterized by the identity

$$\left(x+\sqrt{x^2+s}\right)^n = T_n(x,s)+U_{n-1}(x,s)\sqrt{x^2+s} \tag{1.15}$$

if we set $U_{-1}(x,s)=0.$

The classical Chebyshev polynomials are solutions of the following *eigenvalue problems:*

$$(x^2-1)T_n''(x)+xT_n'(x)=n^2T_n(x) \tag{1.16}$$

and

$$(x^2-1)U_n''(x)+3xU_n'(x)=n(n+2)U_n(x). \tag{1.17}$$

They also satisfy the Rodrigues-type formulae

$$T_n(x)=\frac{(-1)^n}{(2n-1)!!}\sqrt{1-x^2}\left(\frac{d}{dx}\right)^n\left((1-x^2)^{n-\frac{1}{2}}\right) \tag{1.18}$$

and

$$U_n(x)=\frac{(-1)^n(n+1)}{(2n+1)!!}\frac{1}{\sqrt{1-x^2}}\left(\frac{d}{dx}\right)^n\left((1-x^2)^{n+\frac{1}{2}}\right). \tag{1.19}$$

More information and references about the classical Chebyshev polynomials can be found in [7], Section 9.8.2.

## 2. (q,b) - Fibonacci polynomials

We employ the usual abbreviations of $q-$analysis such as the $q-$Pochhammer symbols

$(a;q)_n=(1-a)(1-qa)\cdots(1-q^{n-1}a),$ $(a;q)_\infty=\prod_{j=0}^{\infty}(1-q^ja)$ and

$(a;q)_{-n}=\frac{(a;q)_\infty}{(q^{-n}a;q)_\infty}=\frac{1}{(q^{-n}a;q)_n}.$ The $q-$binomial coefficients are denoted by

$\begin{bmatrix} n \\ k \end{bmatrix}=\begin{bmatrix} n \\ k \end{bmatrix}_q=\frac{(q;q)_n}{(q;q)_k(q;q)_{n-k}}.$ In place of $\begin{bmatrix} n \\ 1 \end{bmatrix}$ we write $[n].$

We will also need the $q-$binomial theorem in the form

$$\sum_{k\geq 0}\frac{(a;q)_k}{(q;q)_k}z^k=\frac{(az;q)_\infty}{(z;q)_\infty}. \tag{2.1}$$

As is well known (cf. e.g. [2] or [5]) the Carlitz $q-$Fibonacci polynomials

$$F_n(x,s,q)=\sum_{k=0}^{\left\lfloor\frac{n-1}{2}\right\rfloor}q^{k^2}\begin{bmatrix}n-1-k\\k\end{bmatrix}s^k x^{n-1-2k} \tag{2.2}$$

satisfy the recursion

$$F_n(x,s,q)=xF_{n-1}(x,qs,q)+qsF_{n-2}(x,q^2s,q) \tag{2.3}$$

with initial values $F_0(x,s,q)=0$ and $F_1(x,s,q)=1.$

Our first aim is the study of the following $q-$analogue of $F_n\left(x,\frac{s}{(1-b)^2}\right).$

**Definition 2.1**

*The polynomials*

$$F_n(x,b,s,q)=\sum_{k=0}^{\left\lfloor\frac{n-1}{2}\right\rfloor}q^{k^2}\begin{bmatrix}n-k-1\\k\end{bmatrix}\frac{1}{(qb;q)_k\left(q^{n-k}b;q\right)_k}s^k x^{n-1-2k} \tag{2.4}$$

*will be called* $(q,b)-$*Fibonacci polynomials.*

These are polynomials in $x$ and $s$ whose coefficients are rational functions in $b.$ It is clear that for $b=0$ they reduce to $F_n(x,0,s,q)=F_n(x,s,q)$.

**Remark**

These polynomials seem to be remarkable since variants of them have occurred in the literature in different contexts. They are closely related to the *Al-Salam and Ismail polynomials* $u_n(x;a,b)$ which originated in [1] and are defined by the recurrence

$$u_n(x;a,b)=x(1+q^{n-1}a)u_{n-1}(x;a,b)-q^{n-2}bu_{n-2}(x;a,b) \tag{2.5}$$

with initial values $u_0(x;a,b)=1$ and $u_1(x;a,b)=(1+a)x.$

It turns out that

$$F_{n+1}(x,b,s,q)=\frac{u_n(x;-qb,-qs)}{(qb;q)_n}. \tag{2.6}$$

Some properties of these polynomials have been studied in [6]. The special case $b=s$ occurred in [2]. M. Schlosser (personal communication) has observed that these polynomials can also be obtained from a specialization of his elliptic binomial Theorem (cf. [8]).

As observed in [1] we have

**Theorem 2.1**

*The* $(q,b)-$*Fibonacci polynomials satisfy the recursion*

$$F_n(x,b,s,q)=xF_{n-1}(x,b,s,q)+\frac{q^{n-2}s}{(1-q^{n-2}b)(1-q^{n-1}b)}F_{n-2}(x,b,s,q) \tag{2.7}$$

*with initial values* $F_0(x,b,s,q)=0$ *and* $F_1(x,b,s,q)=1.$

A closely related recursion gives

**Theorem 2.2**

*The polynomials* $F_n(x,b,s,q)$ *also satisfy the recursion*

$$F_n(x,b,s,q)=xF_{n-1}(x,qb,qs,q)+\frac{qs}{(1-qb)(1-q^2b)}F_{n-2}(x,q^2b,q^2s,q) \tag{2.8}$$

*with initial values* $F_0(x,b,s,q)=0$ *and* $F_1(x,b,s,q)=1.$

Both theorems can be easily verified by comparing coefficients in (2.4).

In order to give more insight into the situation and to indicate the connections with [2] and [8] we will give another proof.

Consider the linear operator $\eta$ defined by

$$\eta f(x,b,s)=f(x,qb,qs) \tag{2.9}$$

for $f\in\mathbb{C}(b,q)[x,s]$.

Let

$$X = x\eta \tag{2.10}$$

and

$$Y = \frac{qs}{(1-qb)(1-q^2b)}\eta^2 \tag{2.11}$$

and let $C_k^n(X,Y)$ be the sum of all words with $k$ letters $Y$ and $n-k$ letters $X$.

Then obviously

$$(X+Y)^n = \sum_{k=0}^{n} C_k^n(X,Y). \tag{2.12}$$

It is easily verified that

$$XY = \frac{1-qb}{1-q^3b} qYX, \tag{2.13}$$

$$Xb = qbX, \tag{2.14}$$

$$Yb = q^2bY. \tag{2.15}$$

As shown in [8] properties (2.13) - (2.15) imply the following $(q,b)$ – binomial theorem:

**Lemma 2.1 (M. Schlosser [8])**

$$(X+Y)^n = \sum_{k=0}^{n} \frac{\left(q^{k+1}b;q\right)_{n-k}\left(q^{k+1};q\right)_{n-k}}{\left(q^{2k+1}b;q\right)_{n-k}\left(q;q\right)_{n-k}} Y^k X^{n-k} = \sum_{k=0}^{n} \begin{bmatrix} n \\ k \end{bmatrix} \frac{\left(q^{k+1}b;q\right)_k}{\left(q^{n+1}b;q\right)_k} Y^k X^{n-k} \tag{2.16}$$

*or equivalently*

$$C_k^n\left(X,Y\right) = \begin{bmatrix} n \\ k \end{bmatrix} \frac{\left(q^{k+1}b;q\right)_k}{\left(q^{n+1}b;q\right)_k} Y^k X^{n-k} = \begin{bmatrix} n \\ k \end{bmatrix} \frac{q^{k^2}}{\left(q^{n+1}b;q\right)_k \left(qb;q\right)_k} s^k x^{n-k} \eta^{n+k}. \tag{2.17}$$

**Proof**

We have to show that $(X+Y)^n=\sum_{k=0}^{n}c(n,k,b)Y^kX^{n-k}$ with

$$c(n,k,b)=\begin{bmatrix} n \\ k \end{bmatrix}\frac{\left(q^{k+1}b;q\right)_k}{\left(q^{n+1}b;q\right)_k}.$$

Since $XY^k=\frac{1-qb}{1-q^{2k+1}b}q^kY^kX$ and (2.14) and (2.15) imply

$$c(n,k,b)=c(n-1,k-1,q^2b)+q^k\frac{1-qb}{1-q^{2k+1}b}c(n-1,k,qb) \tag{2.18}$$

it suffices to verify that

$$\begin{bmatrix} n \\ k \end{bmatrix}\frac{\left(q^{k+1}b;q\right)_k}{\left(q^{n+1}b;q\right)_k}=\begin{bmatrix} n-1 \\ k-1 \end{bmatrix}\frac{\left(q^{k+2}b;q\right)_{k-1}}{\left(q^{n+2}b;q\right)_{k-1}}+\begin{bmatrix} n-1 \\ k \end{bmatrix}\frac{\left(q^{k+2}b;q\right)_k}{\left(q^{n+1}b;q\right)_k}\frac{q^k(1-qb)}{1-q^{2k+1}b}$$

or equivalently

$$\begin{bmatrix} n \\ k \end{bmatrix}\left(1-q^{k+1}b\right)=\begin{bmatrix} n-1 \\ k-1 \end{bmatrix}\left(1-q^{n+1}b\right)+\begin{bmatrix} n-1 \\ k \end{bmatrix}q^k(1-qb).$$

But this is a trivial consequence of the recursions

$$\begin{bmatrix} n \\ k \end{bmatrix}=q^k\begin{bmatrix} n-1 \\ k \end{bmatrix}+\begin{bmatrix} n-1 \\ k-1 \end{bmatrix} \tag{2.19}$$

and

$$\begin{bmatrix} n \\ k \end{bmatrix}=\begin{bmatrix} n-1 \\ k \end{bmatrix}+q^{n-k}\begin{bmatrix} n-1 \\ k-1 \end{bmatrix}. \tag{2.20}$$

for the $q-$binomial coefficients.

If we apply (2.17) to the constant polynomial 1 we get

$$C_k^n(X,Y)1=\begin{bmatrix} n \\ k \end{bmatrix}\frac{\left(q^{k+1}b;q\right)_k}{\left(q^{n+1}b;q\right)_k}Y^kX^{n-k}1=\begin{bmatrix} n \\ k \end{bmatrix}\frac{q^{k^2}}{\left(q^{n+1}b;q\right)_k\left(qb;q\right)_k}s^kx^{n-k}. \tag{2.21}$$

**Proof of Theorem 2.2**

We must show that (2.8) implies (2.4).

Define the length $\lambda\left(Z_1 Z_2 \cdots Z_r\right)$ of a word $Z_1 Z_2 \cdots Z_r$ with $Z_i \in \{X, Y\}$ by $\lambda\left(Z_1 Z_2 \cdots Z_r\right) = \lambda(Z_1) + \cdots + \lambda(Z_r),$ where $\lambda(X) = 1$ and $\lambda(Y) = 2.$

As in [2] define Fibonacci words $\mathbf{F}_n(X,Y)$ as the sum of all words of length $n-1$. We set $\mathbf{F}_0(X,Y) = 0$ and let $\mathbf{F}_1(X,Y) = \varepsilon$ be the empty word so that $\left(\mathbf{F}_n(X,Y)\right)_{n\geq 0} = \left(0, \varepsilon, X, X^2 + Y, X^3 + XY + YX, \cdots\right).$

Then obviously

$$\mathbf{F}_n(X,Y) = X\mathbf{F}_{n-1}(X,Y) + Y\mathbf{F}_{n-2}(X,Y) = \mathbf{F}_{n-1}(X,Y)X + \mathbf{F}_{n-2}(X,Y)Y \tag{2.22}$$

and

$$\mathbf{F}_n(X,Y) = \sum_{k=0}^{\left\lfloor \frac{n-1}{2} \right\rfloor} C_k^{n-1-k}\left(X,Y\right). \tag{2.23}$$

For example

$$\mathbf{F}_5(X,Y) = C_0^4\left(X,Y\right) + C_1^3\left(X,Y\right) + C_2^2\left(X,Y\right) = XXXX + \left(XXY + XYX + YXX\right) + \left(YY\right) = X^4 + X^2Y + XYX + YX^2 + Y^2.$$

Recursion (2.8) means that $F_n(x,b,s,q) = XF_{n-1}(x,b,s,q) + YF_{n-2}(x,b,s,q),$ i.e. by (2.22) and (2.23) that

$$F_n(x,b,s,q) = \sum_{k=0}^{\left\lfloor \frac{n-1}{2} \right\rfloor} C_k^{n-1-k}\left(X,Y\right)1. \tag{2.24}$$

By (2.23) and (2.17) we get

$$\mathbf{F}_n(X,Y) = \sum_{k=0}^{\left\lfloor \frac{n-1}{2} \right\rfloor} C_k^{n-1-k}\left(X,Y\right) = \sum_{k=0}^{\left\lfloor \frac{n-1}{2} \right\rfloor} \begin{bmatrix} n-1-k \\ k \end{bmatrix} \frac{q^{k^2}}{\left(q^{n-k}b;q\right)_k \left(qb;q\right)_k} s^k x^{n-1-2k} \eta^{n-1} \tag{2.25}$$

and therefore

$$F_n(x,b,s,q) = \sum_{k=0}^{\left\lfloor \frac{n-1}{2} \right\rfloor} q^{k^2} \begin{bmatrix} n-k-1 \\ k \end{bmatrix} \frac{1}{\left(qb;q\right)_k \left(q^{n-k}b;q\right)_k} s^k x^{n-1-2k}. \tag{2.26}$$

Thus Theorem 2.2 is proved.

**Proof of Theorem 2.1**

This follows by (2.22) which gives

$$F_n\left(x,b,s,q\right)=\mathbf{F}_n(X,Y)1=\mathbf{F}_{n-1}(X,Y)X+\mathbf{F}_{n-2}(X,Y)Y$$
$$=F_{n-1}\left(x,b,s,q\right)\eta^{n-2}x\eta 1+F_{n-2}\left(x,b,s,q\right)\eta^{n-3}\frac{qs}{(1-qb)(1-q^2b)}\eta^2 1$$
$$=xF_{n-1}\left(x,b,s,q\right)+\frac{q^{n-2}s}{(1-q^{n-2}b)(1-q^{n-1}b)}F_{n-2}\left(x,b,s,q\right).$$

If we extend these polynomials to negative indices by

$$F_{-n}(x,b,s,q)=\frac{(-1)^{n-1}}{s^n}q^{\binom{n+1}{2}}\left(\frac{b}{q^{n-1}};q\right)_n\left(\frac{b}{q^n};q\right)_n F_n\left(x,\frac{b}{q^n},\frac{s}{q^n},q\right)$$

then both recurrences (2.7) and (2.8) remain true.

**Definition 2.2**

*The* $(q,-1)-$*Fibonacci polynomials*

$$F_{n+1}(x,-1,s,q)=\sum_{k=0}^{\left\lfloor\frac{n}{2}\right\rfloor}q^{k^2}\begin{bmatrix}n-k\\k\end{bmatrix}\frac{1}{\left(-q;q\right)_k\left(-q^{n+1-k};q\right)_k}s^k x^{n-2k} \tag{2.27}$$

*or in hypergeometric form*

$$F_{n+1}(x,-1,s,q)=\sum_{k=0}^{\left\lfloor\frac{n}{2}\right\rfloor}\frac{\left(q^{-n};q^2\right)_k\left(q^{1-n};q^2\right)_k}{\left(q^{-2n};q^2\right)_k\left(q^2;q^2\right)_k}\left(-s\right)^k x^{n-2k} \tag{2.28}$$

*are called generalized* $q-$*Fibonacci polynomials.*

We now consider analogues of the Fibonacci matrices. Let

$$C(x,b,s,q)=\begin{pmatrix} 0 & 1 \\ \dfrac{s}{(1-b)(1-qb)} & x \end{pmatrix}. \tag{2.29}$$

Then we have

$$\begin{pmatrix} \dfrac{s}{(1-qb)(1-q^2b)}F_{n-1}(x,q^2b,q^2s,q) & F_n(x,qb,qs,q) \\ \dfrac{s}{(1-qb)(1-q^2b)}F_n(x,q^2b,q^2s,q) & F_{n+1}(x,qb,qs,q) \end{pmatrix}\begin{pmatrix} 0 & 1 \\ \dfrac{s}{(1-b)(1-qb)} & x \end{pmatrix}$$

$$=\begin{pmatrix} \dfrac{s}{(1-b)(1-qb)}F_n(x,qb,qs,q) & F_{n+1}(x,b,s,q) \\ \dfrac{s}{(1-b)(1-qb)}F_{n+1}(x,qb,qs,q) & F_{n+2}(x,b,s,q) \end{pmatrix}$$

and therefore we get the $(q,b)-$ *Fibonacci matrices*

$$C\left(x,q^{n-1}b,q^{n-1}s,q\right)\cdots C\left(x,b,s,q\right)=\begin{pmatrix} \dfrac{s}{(1-b)(1-qb)}F_{n-1}(x,qb,qs,q) & F_n(x,b,s,q) \\ \dfrac{s}{(1-b)(1-qb)}F_n(x,qb,qs,q) & F_{n+1}(x,b,s,q) \end{pmatrix}. \tag{2.30}$$

Taking determinants we get the following $(q,b)-$ analogue of Cassini's formula.

**Theorem 2.3 ( (q,b) -Cassini formula)**

*For all* $n\in\mathbb{Z}$

$$F_{n-1}(x,qb,qs,q)F_{n+1}(x,b,s,q)-F_n(x,b,s,q)F_n(x,qb,qs,q)=(-1)^n\frac{q^{\binom{n}{2}}}{\left(qb;q\right)_{n-1}\left(q^2b;q\right)_{n-1}}s^{n-1}. \tag{2.31}$$

In fact a slightly more general theorem is true.

**Theorem 2.4 ( (q,b) -Cassini-Euler formula,[6])**

*Let*

$$d(n,k,b,s)=\frac{s}{(1-b)(1-qb)}\left(F_{n-1}(x,qb,qs,q)F_{n+k}(x,b,s,q)-F_{n+k-1}(x,qb,qs,q)F_{n}(x,b,s,q)\right). \quad (2.32)$$

*Then*

$$d(n,k,b,s)=q^{\binom{n}{2}}\frac{(-s)^n}{(b;q)_n(qb;q)_n}F_k\left(x,q^n b,q^n s,q\right). \quad (2.33)$$

**Proof**

This has been proved with other methods in [6], Theorem 3.1.

In our approach we have only to verify that

$$\begin{pmatrix} \frac{s}{(1-qb)(1-q^2b)}F_{n-1}(x,q^2b,q^2s,q) & F_n(x,qb,qs,q) \\ \frac{s}{(1-qb)(1-q^2b)}F_{n+k-1}(x,q^2b,q^2s,q) & F_{n+k}(x,qb,qs,q) \end{pmatrix}\begin{pmatrix} 0 & 1 \\ \frac{s}{(1-b)(1-qb)} & x \end{pmatrix}$$
$$=\begin{pmatrix} \frac{s}{(1-b)(1-qb)}F_{n}(x,qb,qs,q) & F_{n+1}(x,b,s,q) \\ \frac{s}{(1-b)(1-qb)}F_{n+k}(x,qb,qs,q) & F_{n+k+1}(x,b,s,q) \end{pmatrix}.$$

Taking determinants this gives

$$d(n,k,b,s)=\frac{-s}{(1-b)(1-qb)}d(n-1,k,qb,qs).$$

Since

$$d(1,k,b,s)=\frac{-s}{(1-b)(1-qb)}\left(F_k(x,qb,qs,q)\right)$$

the result follows.

## 3. (q,b) -Lucas polynomials

For $b=0$ the Carlitz $q-$Lucas polynomials coincide with the trace of the Fibonacci matrices. Thus a possible choice is to define a $(q,b)-$ analogue of the Lucas polynomials in the same way, i.e.

$$l_n(x,b,s,q)=tr\left(C\left(x,q^{n-1}b,q^{n-1}s,q\right)\cdots C\left(x,b,s,q\right)\right) \tag{3.1}$$

or

$$l_n(x,b,s,q)=F_{n+1}(x,b,s,q)+\frac{s}{(1-b)(1-qb)}F_{n-1}(x,qb,qs,q). \tag{3.2}$$

It is then easily verified that for $n>0$

$$l_n(x,b,s,q)=\sum_{k=0}^{\left\lfloor \frac{n}{2}\right\rfloor} q^{k^2-k}s^k x^{n-2k}\frac{[n]}{[n-k]}\begin{bmatrix} n-k \\ k \end{bmatrix}\frac{1}{(b;q)_k\left(q^{n-k+1}b;q\right)_k}. \tag{3.3}$$

Note that we have $l_0(x,b,s,q)=2$ since $F_{-1}(x,qb,qs)=\dfrac{(1-b)(1-qb)}{s}.$

For $n<0$ we get

$$l_{-n}(x,b,s,q)=(-1)^n\frac{q^{\binom{n+1}{2}}}{s^n}\left(\frac{b}{q^n};q\right)_n\left(\frac{b}{q^{n-1}};q\right)_n l_n\left(x,\frac{b}{q^n},\frac{s}{q^n},q\right). \tag{3.4}$$

Unfortunately there is no linear functional with respect to which they are orthogonal.

For consider the $q-$Lucas polynomials $l_n(x,q)=l_n(x,0,s,q)$ the first terms of which are $2,x,x^2+[2]s,x^3+[3]sx,x^4+[4]sx^2+(q^2+q^4)s^2.$

Here we get $l_1(x,q)l_3(x,q)=l_4(x,q)-q^3sl_2(x,q)-q^2(1-q)s^2.$ For each linear functional $\Lambda$ with $\Lambda\left(l_n(x,s)\right)=0$ for $n>0$ we get $\Lambda\left(l_1(x,q)l_3(x,q)\right)=-q^2(1-q)s^2\neq 0$ if $q\neq 1.$

Fortunately there is another class of polynomials with both a beautiful formula and a 3-term recursion which is in close connection with the $(q,b)-$ Fibonacci polynomials.

**Definition 3.1**

*The* $(q,b)-$*Lucas polynomials* $L_n(x,b,s,q)$ *are defined by*

$$L_n(x,b,s,q)=xL_{n-1}(x,qb,qs,q)+\frac{qs}{(1-qb)(1-q^2b)}L_{n-2}(x,q^2b,q^2s,q) \tag{3.5}$$

*with initial values* $L_0(x,b,s,q)=1-b$ *and* $L_1(x,b,s,q)=x.$

These polynomials are related to the $(q,b)-$ Fibonacci polynomials by

**Theorem 3.1**

$$L_n(x,b,s,q)=F_{n+1}(x,b,s,q)-\frac{q^{2n-1}sb}{\left(1-q^{n-1}b\right)\left(1-q^nb\right)}F_{n-1}(x,b,s,q). \tag{3.6}$$

For this is true for $n=1$ and $n=2.$ Since $L_n(x,b,s,q)$ and $F_{n+1}(x,b,s,q)$ satisfy the same recurrence (3.5) it suffices to show that $\frac{q^{2n-1}sb}{\left(1-q^{n-1}b\right)\left(1-q^nb\right)}F_{n-1}(x,b,s,q)$ satisfies this recursion too. But this obvious from

$$\frac{q^{2n-1}sb}{\left(1-q^{n-1}b\right)\left(1-q^nb\right)}F_{n-1}(x,b,s,q)=x\frac{q^{2n-1}sb}{\left(1-q^{n-1}b\right)\left(1-q^nb\right)}F_{n-2}(x,qb,qs,q)$$
$$+\frac{qs}{(1-qb)(1-q^2b)}\frac{q^{2n-1}sb}{\left(1-q^{n-1}b\right)\left(1-q^nb\right)}F_{n-3}(x,q^2b,q^2s,q).$$

This implies

$$L_n(x,b,s,q)=\sum_{k=0}^{\left\lfloor\frac{n}{2}\right\rfloor}q^{k^2}s^kx^{n-2k}\frac{\begin{bmatrix}n-k\\k\end{bmatrix}-q^{n-k}b\begin{bmatrix}n-1-k\\k-1\end{bmatrix}}{(qb;q)_k\left(q^{n-k}b;q\right)_k}. \tag{3.7}$$

For the coefficient of $s^kx^{n-2k}$ in (3.6) is

$$q^{k^2}\begin{bmatrix} n-k \\ k \end{bmatrix}\frac{1}{(qb;q)_k\left(q^{n+1-k}b;q\right)_k}-\frac{q^{2n-1}b}{\left(1-q^{n-1}b\right)\left(1-q^nb\right)}q^{(k-1)^2}\begin{bmatrix} n-1-k \\ k-1 \end{bmatrix}\frac{1}{(qb;q)_{k-1}\left(q^{n-k}b;q\right)_{k-1}}$$

$$=q^{k^2}\begin{bmatrix} n-k \\ k \end{bmatrix}\frac{1}{(qb;q)_k\left(q^{n-k}b;q\right)_k}\left(\frac{1-q^{n-k}b}{1-q^nb}-q^{2n-2k}b\frac{(1-q^kb)}{\left(1-q^nb\right)}\frac{1-q^k}{1-q^{n-k}}\right)$$

$$=q^{k^2}\frac{1}{(qb;q)_k\left(q^{n-k}b;q\right)_k}\left(\begin{bmatrix} n-k \\ k \end{bmatrix}-q^{n-k}b\begin{bmatrix} n-1-k \\ k-1 \end{bmatrix}\right).$$

**Theorem 3.2**

*The* $(q,b)-$*Lucas polynomials also satisfy the recursion*

$$L_n(x,b,s,q)=xL_{n-1}(x,b,s,q)+\frac{q^{n-1}s}{(1-q^{n-2}b)(1-q^{n-1}b)}L_{n-2}(x,b,s,q). \tag{3.8}$$

This follows from the fact that each term $h_n(x,b,s,q)$ satisfies

$$h_n(x,b,s,q)=xh_{n-1}(x,qb,qs,q)+\frac{qs}{(1-qb)(1-q^2b)}h_{n-2}(x,q^2b,q^2s,q).$$

## 4. Generalized q-Lucas polynomials

The $(q,b)-$Lucas polynomials have very interesting properties if $b=-1.$

For $b=-1$ formula (3.7) reduces to

$$L_n(x,-1,s,q)=F_{n+1}(x,-1,s,q)+\frac{q^{2n-1}s}{\left(1+q^{n-1}\right)\left(1+q^n\right)}F_{n-1}(x,-1,s,q). \tag{4.1}$$

In this case (3.7) reduces to the following $q-$analogue of $L_n\left(x,\frac{s}{4}\right).$

**Theorem 4.1**

*For* $n > 0$

$$L_n(x,-1,s,q) = \sum_{k=0}^{\left\lfloor \frac{n}{2} \right\rfloor} q^{k^2} s^k x^{n-2k} \frac{[n]}{[n-k]} \begin{bmatrix} n-k \\ k \end{bmatrix} \frac{1}{(-q;q)_k \left(-q^{n-k};q\right)_k}. \tag{4.2}$$

**Definition 4.1**

*The polynomials* $L_n(x,-1,s,q)$ *are called generalized* $q-$*Lucas polynomials.*

In hypergeometric form formula (4.2) can be written as

$$L_n(x,-1,s,q) = \sum_{k=0}^{\left\lfloor \frac{n}{2} \right\rfloor} \frac{\left(q^{-n};q^2\right)_k \left(q^{1-n};q^2\right)_k}{\left(q^{-2n+2};q^2\right)_k \left(q^2;q^2\right)_k} \left(-q^2 s\right)^k x^{n-2k}. \tag{4.3}$$

An equivalent formula to (4.1) is

$$L_n(x,-1,s,q) = (1+q^n)F_{n+1}(x,-1,s,q) - q^n x F_n(x,-1,s,q). \tag{4.4}$$

For the right-hand side is $F_{n+1}(x,-1,s,q) + q^n \left(F_{n+1}(x,-1,s,q) - xF_n(x,-1,s,q)\right)$ which by (2.7) coincides with (4.1).

Comparing coefficients we also get

$$q^n L_n(x,-1,s,q) = (1+q^n)F_{n+1}(x,-1,q^2 s,q) - xF_n(x,-1,q^2 s,q). \tag{4.5}$$

Again these polynomials can be extended to negative indices such that the recurrences (3.8) and (3.5) remain true:

$$L_{-n}(x,-1,s,q) = (-1)^n \frac{1}{q^{\binom{n+1}{2}} s^n} (-q;q)_n (-1;q)_n L_n(x,-1,s,q). \tag{4.6}$$

From (2.7) and (3.8) follows that both the generalized $q-$Fibonacci polynomials and the generalized $q-$Lucas polynomials are orthogonal with respect to some linear functional. We want to compute their moments.

For $q=1$ the Fibonacci polynomials $F_n(x,s)$ satisfy

$$x^n = \sum_{k=0}^{\left\lfloor \frac{n}{2} \right\rfloor} \left( \binom{n}{k} - \binom{n}{k-1} \right) (-s)^k F_{n+1-2k}(x,s). \tag{4.7}$$

Therefore the corresponding moments are given by

$$\Lambda\left(x^{2n}\right) = (-s)^n \left( \binom{2n}{n} - \binom{2n}{n-1} \right) = (-s)^n \frac{1}{n+1}\binom{2n}{n}, \quad \Lambda\left(x^{2n+1}\right) = 0. \tag{4.8}$$

For $s=-1$ the even moments are the Catalan numbers $\Lambda\left(x^{2n}\right) = C_n = \frac{1}{n+1}\binom{2n}{n}.$

The corresponding moments for $F_n(x,s,q)$ are $\Lambda\left(x^{2n}\right) = (-qs)^n C_n(q)$, where the (Carlitz-) $q-$ Catalan numbers $C_n(q)$ are defined by $C_n(q) = \sum_{k=0}^{n-1} q^k C_k(q) C_{n-1-k}(q)$ with $C_0(q)=1.$ This result seems to be “folklore”. A proof may be found in [4].

Unfortunately there is no closed formula for these $q-$Catalan numbers. As Wadim Zudilin [9] has shown they are not even $q-$holonomic.

Therefore for general $b$ there can be no simple analogue of (4.7) for $F_n(x,b,s,q).$ But for $b=-1$ such a $q-$analogue exists.

Here we get

**Theorem 4.2**

$$x^n = \sum_{k=0}^{\left\lfloor \frac{n}{2} \right\rfloor} \left( \begin{bmatrix} n \\ k \end{bmatrix} - \begin{bmatrix} n \\ k-1 \end{bmatrix} \right) \frac{(-s)^k}{(-q;q)_k \left(-q^{n+2-2k};q\right)_k} F_{n+1-2k}(x,-1,s,q). \tag{4.9}$$

**Proof**

We verify this by induction. For $n=0$ it is trivially true.

Since $xF_{n+1-2k}(x,-1,s,q) = F_{n+2-2k}(x,-1,s,q) - \frac{q^{n-2k}s}{(1+q^{n-2k})(1+q^{n+1-2k})} F_{n-2k}(x,-1,s,q)$

we have to show that

$$\sum_{k=0}^{\left\lfloor\frac{n+1}{2}\right\rfloor}\left(\begin{bmatrix} n+1 \\ k \end{bmatrix}-\begin{bmatrix} n+1 \\ k-1 \end{bmatrix}\right)\frac{(-s)^k}{(-q;q)_k\left(-q^{n+3-2k};q\right)_k}F_{n+2-2k}(x,-1,s,q)$$

$$=\sum_{k=0}^{\left\lfloor\frac{n}{2}\right\rfloor}\left(\begin{bmatrix} n \\ k \end{bmatrix}-\begin{bmatrix} n \\ k-1 \end{bmatrix}\right)\frac{(-s)^k}{(-q;q)_k\left(-q^{n+2-2k};q\right)_k}F_{n+2-2k}(x,-1,s,q)$$

$$-\sum_{k=0}^{\left\lfloor\frac{n}{2}\right\rfloor}\left(\begin{bmatrix} n \\ k \end{bmatrix}-\begin{bmatrix} n \\ k-1 \end{bmatrix}\right)\frac{(-s)^k}{(-q;q)_k\left(-q^{n+2-2k};q\right)_k}\frac{q^{n-2k}s}{(1+q^{n-2k})(1+q^{n+1-2k})}F_{n-2k}(x,-1,s,q)$$

This is equivalent with

$$\left(1+q^{n+2-2k}\right)\left(\begin{bmatrix} n+1 \\ k \end{bmatrix}-\begin{bmatrix} n+1 \\ k-1 \end{bmatrix}\right)=\left(1+q^{n+2-k}\right)\left(\begin{bmatrix} n \\ k \end{bmatrix}-\begin{bmatrix} n \\ k-1 \end{bmatrix}\right)+(1+q^k)q^{n+2-2k}\left(\begin{bmatrix} n \\ k-1 \end{bmatrix}-\begin{bmatrix} n \\ k-2 \end{bmatrix}\right),$$

which is easily verified.

Let $\Lambda_F$ be the linear functional corresponding to the sequence $\left(F_{n+1}(x,-1,s,q)\right)_{n\geq 0}$.

Then

$$\Lambda_F(X^{2n})=\frac{1}{[n+1]}\begin{bmatrix} 2n \\ n \end{bmatrix}\frac{(-s)^n}{(1+q)(1+q^{n+1})\prod_{j=2}^{n}(1+q^j)^2}=\frac{1}{[n+1]}\begin{bmatrix} 2n \\ n \end{bmatrix}\frac{(-s)^n}{(-q;q)_n(-q^2;q)_n}$$
$$=(-s)^n q^{n^2}(1+q)\begin{bmatrix} \frac{1}{2} \\ n+1 \end{bmatrix}_{q^2}. \tag{4.10}$$

and $\Lambda_F(X^{2n+1})=0.$

There is a corresponding identity for the generalized $q-$Lucas polynomials.

**Theorem 4.3**

*Let* $L_n^*(x,s,q)=L_n(x,-1,s,q)$ *for* $n>0$ *and* $L_0^*(x,s,q)=1.$

*Then*

$$\sum_{k=0}^{\left\lfloor\frac{n}{2}\right\rfloor}\begin{bmatrix} n \\ k \end{bmatrix}\frac{(-qs)^k}{(-q;q)_k\left(-q^{n-2k+1};q\right)_k}L_{n-2k}^*(x,s,q)=x^n. \tag{4.11}$$

With other words

$$\sum_{k=0}^{n}\begin{bmatrix}2n\\k\end{bmatrix}\frac{(-qs)^k}{(-q;q)_k\left(-q^{2n-2k+1};q\right)_k}L_{2n-2k}(x,-1,s,q)=x^{2n}+\begin{bmatrix}2n\\n\end{bmatrix}\frac{(-qs)^n}{(-q;q)_n\left(-q;q\right)_n} \tag{4.12}$$

and

$$\sum_{k=0}^{n}\begin{bmatrix}2n+1\\k\end{bmatrix}\frac{(-qs)^k}{(-q;q)_k\left(-q^{2n-2k+2};q\right)_k}L_{2n+1-2k}(x,-1,s,q)=x^{2n+1}. \tag{4.13}$$

Both identities are trivial for $n=0.$ Now suppose they have been proved for $n-1.$

From (3.8) we see that

$$L_{n+1-2k}(x,-1,s,q)=xL_{n-2k}(x,-1,s,q)+\frac{q^{n-2k}s}{(1+q^{n-1-2k})(1+q^{n-2k})}L_{n-1-2k}(x,-1,s,q)$$

Therefore (4.13) implies $\sum_{k=0}^{n}\begin{bmatrix}2n-1\\k\end{bmatrix}\frac{(-qs)^k}{(-q;q)_k\left(-q^{2n-2k};q\right)_k}xL_{2n-1-2k}(x,-1,s,q)=x^{2n}.$

Thus (4.12) is equivalent with

$$\sum_{k=0}^{n}\begin{bmatrix}2n\\k\end{bmatrix}\frac{(-qs)^k}{(-q;q)_k\left(-q^{2n-2k+1};q\right)_k}L_{2n-2k}(x,-1,s,q)$$

$$=\begin{bmatrix}2n\\n\end{bmatrix}\frac{(-qs)^n}{(-q;q)_n\left(-q;q\right)_n}+\sum_{k=0}^{n}\begin{bmatrix}2n-1\\k\end{bmatrix}\frac{(-qs)^k}{(-q;q)_k\left(-q^{2n-2k};q\right)_k}xL_{2n-1-2k}(x,-1,s,q)$$

$$=\begin{bmatrix}2n\\n\end{bmatrix}\frac{(-qs)^n}{(-q;q)_n\left(-q;q\right)_n}+\sum_{k=0}^{n}\begin{bmatrix}2n-1\\k\end{bmatrix}\frac{(-qs)^k}{(-q;q)_k\left(-q^{2n-2k};q\right)_k}L_{2n-2k}(x,-1,s,q)$$

$$-\sum_{k=0}^{n}\begin{bmatrix}2n-1\\k\end{bmatrix}\frac{(-qs)^k}{(-q;q)_k\left(-q^{2n-2k};q\right)_k}\frac{q^{2n-2k-1}s}{(1+q^{2n-2-2k})(1+q^{2n-1-2k})}L_{2n-2-2k}(x,-1,s,q)$$

$$=\begin{bmatrix}2n\\n\end{bmatrix}\frac{(-qs)^n}{(-q;q)_n\left(-q;q\right)_n}+\sum_{k=0}^{n}\begin{bmatrix}2n-1\\k\end{bmatrix}\frac{(-qs)^k}{(-q;q)_k\left(-q^{2n-2k};q\right)_k}L_{2n-2k}(x,-1,s,q)$$

$$+\sum_{k=0}^{n}\begin{bmatrix}2n-1\\k-1\end{bmatrix}\frac{(-qs)^k}{(-q;q)_{k-1}\left(-q^{2n-2k+2};q\right)_{k-1}}\frac{q^{2n-2k}}{(1+q^{2n-2k})(1+q^{2n+1-2k})}L_{2n-2k}(x,-1,s,q)$$

or

$$\sum_{k=0}^{n}\begin{bmatrix}2n\\k\end{bmatrix}\frac{(-qs)^k}{(-q;q)_k\left(-q^{2n-2k+1};q\right)_k}L_{2n-2k}(x,-1,s,q)$$

$$=\begin{bmatrix}2n\\n\end{bmatrix}\frac{(-qs)^n}{(-q;q)_n(-q;q)_n}+\sum_{k=0}^{n}\begin{bmatrix}2n-1\\k\end{bmatrix}\frac{(-qs)^k}{(-q;q)_k\left(-q^{2n-2k};q\right)_k}L_{2n-2k}(x,-1,s,q)$$

$$+\sum_{k=0}^{n}\begin{bmatrix}2n-1\\k-1\end{bmatrix}\frac{(-qs)^k}{(-q;q)_{k-1}\left(-q^{2n-2k+2};q\right)_{k-1}}\frac{q^{2n-2k}}{(1+q^{2n-2k})(1+q^{2n+1-2k})}L_{2n-2k}(x,-1,s,q)$$

or

$$\begin{bmatrix}2n\\k\end{bmatrix}\left(1+q^{2n-2k}\right)=\begin{bmatrix}2n-1\\k\end{bmatrix}\left(1+q^{2n-k}\right)+\begin{bmatrix}2n-1\\k-1\end{bmatrix}\left(q^{2n-k}+q^{2n-2k}\right)$$

This is true because of the recursions of the $q-$binomial coefficients.

In the same way (4.13) is obtained.

Let $\Lambda_L$ be the linear functional corresponding to the sequence $\left(L_n^*(x,s,q)\right)_{n\geq 0}$.

Then

$$\Lambda_L(x^{2n})=\begin{bmatrix}2n\\n\end{bmatrix}\frac{(-qs)^n}{(-q;q)_n^2}. \tag{4.14}$$

## 5. q-Chebyshev polynomials

Now we are ready to define the announced $q-$analogues of the Chebyshev polynomials.

**Definition 5.1**

*The polynomials*

$$U_n(x,s,q)=F_{n+1}(x,-1,s,q)\left(-q;q\right)_n=\sum_{k=0}^{\left\lfloor\frac{n}{2}\right\rfloor}q^{k^2}\begin{bmatrix}n-k\\k\end{bmatrix}\left(1+q^{k+1}\right)\cdots\left(1+q^{n-k}\right)s^kx^{n-2k} \tag{5.1}$$

*are called $q-$Chebyshev polynomials of the second kind.*

**Theorem 5.1**

*The $q-$Chebyshev polynomials of the second kind satisfy*

$$U_n(x,s,q)=\left(1+q^n\right)xU_{n-1}(x,s,q)+q^{n-1}sU_{n-2}(x,s,q) \tag{5.2}$$

*with initial values $U_0(x,s,q)=1$ and $U_1(x,s,q)=(1+q)x.$*

The proof is an immediate consequence of (2.7).

The first terms of the sequence $\left(U_n(x,s,q)\right)_{n\geq 0}$ are
$1,\ (1+q)x,\ (1+q)(1+q^2)x^2+qs,\ (1+q)(1+q^2)(1+q^3)x^3+q(1+q)(1+q^2)sx,\cdots.$

It is clear that $U_n(x,-1,1)=U_n(x)$ is the classical Chebyshev polynomial of the second kind.

**Remark**

As already mentioned these polynomials are special cases of the Al-Salam and Ismail polynomials. More precisely we have

$$U_n(x,s,q)=u_n(x;q,-qs) \tag{5.3}$$

which can be easily verified by comparing (5.2) with (2.5).

**Definition 5.2**

*The polynomials*

$$T_n(x,s,q)=\left(-q;q\right)_{n-1}L_n(x,-1,s,q)=\sum_{k=0}^{\left\lfloor \frac{n}{2}\right\rfloor}q^{k^2}\frac{[n]}{[n-k]}\begin{bmatrix}n-k\\k\end{bmatrix}\frac{\left(-q;q\right)_{n-1}}{\left(-q;q\right)_k\left(-q^{n-k};q\right)_k}s^kx^{n-2k} \tag{5.4}$$

*for $n>0$ and $T_0(x,s,q)=1$ are called $q-$Chebyshev polynomials of the first kind.*

**Theorem 5.2**

*The $q-$Chebyshev polynomials of the first kind satisfy*

$$T_n(x,s,q)=(1+q^{n-1})xT_{n-1}(x,s,q)+q^{n-1}sT_{n-2}(x,s,q) \tag{5.5}$$

*for $n\geq 2$ with initial values $T_0(x,s,q)=1$ and $T_1(x,s,q)=x.$*

The proof follows from (3.8).

The first terms of the sequence $\left(T_n(x,s,q)\right)_{n\geq 0}$ are

$$1,\ x,\ (1+q)x^2+qs,\ (1+q)(1+q^2)x^3+q(1+q+q^2)sx,$$
$$(1+q)(1+q^2)(1+q^3)x^4+q(1+q)(1+q^2)^2sx^2+q^4s^2,\cdots.$$

The polynomial $T_n(x,-1,1)=T_n(x)$ is the classical Chebyshev polynomial of the first kind.

For negative indices we get

$$U_{-n-2}(x,s,q)=(-1)^n\left(\frac{q}{s}\right)^{n+1}U_n(x,s) \tag{5.6}$$

and

$$T_{-n}(x,s,q)=\frac{(-1)^n}{s^n}T_n(x,s,q). \tag{5.7}$$

Several formulae for the classical Chebyshev polynomials have beautiful $q-$analogues:

**Theorem 5.3**

*The $q-$Chebyshev polynomials satisfy the recursion*

$$\begin{pmatrix} T_{n+1}(x,s,q) \\ U_n(x,s,q) \end{pmatrix} = \begin{pmatrix} q^n x & (x^2+qs)\eta^2 \\ 1 & q^n x \end{pmatrix} \begin{pmatrix} T_n(x,s,q) \\ U_{n-1}(x,s,q) \end{pmatrix} \tag{5.8}$$

*with initial values $U_{-1}(x,s,q)=0$ and $T_0(x,s,q)=1$ and are uniquely determined by this condition.*

To prove (5.8) we must show that

$$T_n(x,s,q)=U_n(x,s,q)-q^nxU_{n-1}(x,s,q) \tag{5.9}$$

and

$$T_{n+1}(x,s,q)=q^nxT_n(x,s,q)+(x^2+qs)U_{n-1}(x,q^2s,q) \tag{5.10}$$

or equivalently

$$(1+q^n)L_{n+1}(x,-1,s,q)-q^n xL_n(x,-1,s,q)=(x^2+qs)F_n(x,-1,q^2s,q). \qquad (5.11)$$

(5.9) is equivalent with (4.4).

(5.11) follows by comparing coefficients in (2.7) and (4.2).

The uniqueness is obvious.

**Theorem 5.4**

*Let* $A=\sqrt{(x^2+qs)\eta^2}$ *be the (symbolic) square-root of the linear operator* $(x^2+qs)\eta^2$ *and let* $p_n(x,A)=(x+A)(qx+A)\cdots(q^{n-1}x+A).$

*Then we get the Binet-like formulae*

$$T_n(x,s,q)=\frac{p_n(x,A)+p_n(x,-A)}{2}1=\sum_{k=0}^{\left\lfloor \frac{n}{2}\right\rfloor} q^{\binom{n-2k}{2}}\begin{bmatrix} n \\ 2k\end{bmatrix} x^{n-2k}\prod_{j=0}^{k-1}\left(x^2+q^{2j+1}s\right) \qquad (5.12)$$

*and*

$$U_n(x,s,q)=\frac{p_{n+1}(x,A)-p_{n+1}(x,-A)}{2A}1=\sum_{k=0}^{\left\lfloor \frac{n+1}{2}\right\rfloor} q^{\binom{n-2k}{2}}\begin{bmatrix} n+1 \\ 2k+1\end{bmatrix} x^{n-2k}\prod_{j=0}^{k-1}\left(x^2+q^{2j+1}s\right). \qquad (5.13)$$

*Equivalently*

$$p_n(x,A)=T_n(x,s,q)+U_{n-1}(x,s,q)A. \qquad (5.14)$$

**Proof**

To prove this we have only to show that (5.8) holds.

But this follows immediately from

$$p_{n+1}(x,A)+p_{n+1}(x,-A)=(q^n x+A)p_n(x,A)+(q^n x-A)p_n(x,-A)$$
$$=q^n x\left(p_n(x,A)+p_n(x,-A)\right)+A^2\,\frac{p_n(x,A)-p_n(x,-A)}{A}$$

and

$$\frac{p_{n+1}(x,A)-p_{n+1}(x,-A)}{A}=\frac{(q^n x+A)p_n(x,A)-(q^n x-A)p_n(x,-A)}{A}$$
$$=q^n x\frac{p_n(x,A)-p_n(x,-A)}{A}+p_n(x,A)+p_n(x,-A).$$

Here we used the $q-$binomial theorem $(x+y)(qx+y)\cdots(q^{n-1}x+y)=\sum_{k=0}^{n} q^{\binom{k}{2}}\begin{bmatrix} n \\ k \end{bmatrix} x^k y^{n-k}.$

**Remark**

Note that for $(s,q)=(-1,1)$ (5.14) reduces to (1.15).

In this case the choice $x=\cos\vartheta$ gives the trigonometric property $T_n(\cos\vartheta)=\cos(n\vartheta)$. Unfortunately this property does not seem to have a simple $q-$analogue.

**Theorem 5.5**

*Let* $V_n(x,s,q)=\begin{pmatrix} q^n x & q^n\left(x^2+qs\right) \\ 1 & x \end{pmatrix}.$

*Then*

$$\begin{pmatrix} T_n(x,s,q) & (x^2+qs)U_{n-1}(x,q^2s,q) \\ U_{n-1}(x,qs,q) & T_n(x,qs,q) \end{pmatrix}=V_0(x,s,q)V_1(x,s,q)\cdots V_{n-1}(x,s,q). \qquad (5.15)$$

**Proof**

It suffices to show that

$$\begin{pmatrix} T_n(x,s,q) & (x^2+qs)U_{n-1}(x,q^2s,q) \\ U_{n-1}(x,qs,q) & T_n(x,qs,q) \end{pmatrix}=\begin{pmatrix} T_{n-1}(x,s,q) & (x^2+qs)U_{n-2}(x,q^2s,q) \\ U_{n-2}(x,qs,q) & T_{n-1}(x,qs,q) \end{pmatrix}\begin{pmatrix} q^{n-1}x & q^{n-1}(x^2+qs) \\ 1 & x \end{pmatrix}.$$

This follows from (5.9), (5.10) and (4.5).

Taking determinants and observing that $T_n(qx,q^2s,q)=q^nT_n(x,s,q)$ and $U_n(qx,q^2s,q)=q^nU_n(x,s,q)$ we get

**Theorem 5.6**

$$T_n(x,s,q)T_n(x,qs,q)-(x^2+qs)U_{n-1}(x,qs,q)U_{n-1}(x,q^2s,q)=q^{\binom{n+1}{2}}(-s)^n \qquad (5.16)$$

*or equivalently*

$$T_n(qx,s,q)T_n\left(\sqrt{q}x,s,q\right)-q^{\frac{2n-1}{2}}(qx^2+s)U_{n-1}(x,s,q)U_{n-1}\left(\sqrt{q}x,s,q\right)=q^{\frac{n^2}{2}}(-s)^n. \qquad (5.17)$$

These are $q-$analogues of the well-known formula

$$T_n^2(x)-(x^2-1)U_{n-1}^2(x)=1.$$

We can also deduce $q-$analogues of the differential equations of the Chebyshev polynomials.

Let $D$ be the $q-$differentiation operator defined by $Df(x)=\frac{f(x)-f(qx)}{x-qx}$.

**Theorem 5.7**

$$DT_n(x,s,q)=[n]U_{n-1}(x,s,q). \qquad (5.18)$$

This is obvious by comparing the formulae for $F_n(x,-1,s,q)$ and $L_n(x,-1,s,q)$.

**Theorem 5.8**

$$(x^2+qs)DU_{n-1}(x,q^2s,q)+q^{n-1}xU_{n-1}(x,s,q)=[n]T_n(x,s,q). \qquad (5.19)$$

**Proof**

The assertion is equivalent with

$$(x^2+qs)DF_n(x,q^2s,q)+q^{n-1}xF_n(x,s,q)=[n]L_n(x,s,q)$$

or

$$(x^2+qs)\sum q^{k^2}\begin{bmatrix} n-k-1 \\ k \end{bmatrix}[n-1-2k]\frac{1}{(-q;q)_k\left(-q^{n-k};q\right)_k}q^{2k}s^k x^{n-2-2k}$$

$$+q^{n-1}\sum q^{k^2}\begin{bmatrix} n-k-1 \\ k \end{bmatrix}\frac{1}{(-q;q)_k\left(-q^{n-k};q\right)_k}s^k x^{n-2k}$$

$$=[n]\sum q^{k^2}\begin{bmatrix} n-k \\ k \end{bmatrix}\frac{[n]}{[n-k]}\frac{1}{(-q;q)_k\left(-q^{n-k};q\right)_k}s^k x^{n-2k}$$

or by comparing coefficients of $q^{k^2}\begin{bmatrix} n-k \\ k \end{bmatrix}\frac{1}{(-q;q)_k\left(-q^{n-k};q\right)_k}s^k x^{n-2k}$

with

$$q^{2k}[n-2k][n-1-2k]+[k][n-k](1+q^k)(1+q^{n-k})+q^{n-1}[n-2k]=[n]^2,$$

which is easily verified.

Combining these Theorems and observing the product formula $D\left(f(x)g(x)\right)=f(x)Dg(x)+g(qx)Df(x)$ and the fact that $U_n(qx,q^2s,q)=q^nU_n(x,s,q)$ we get

**Theorem 5.9**

*The $q-$Chebyshev polynomials satisfy the $q-$differential equations*

$$\left(x^2+qs\right)D^2T_n(x,q^2s,q)+q^{n-1}xDT_n(x,s,q)=[n]^2T_n(x,s,q) \tag{5.20}$$

*and*

$$\left(x^2+qs\right)D^2U_n(x,q^2s,q)+q^{n-1}[3]xDU_n(x,s,q)=[n][n+2]U_n(x,s,q). \tag{5.21}$$

**Remark**

An equivalent formulation of (5.20) is

$$\left(qx^2+\frac{s}{q}\right)D^2T_n(x,s,q)+xDT_n(x,s,q)=q^{-n}[n]^2T_n(qx,s,q). \tag{5.22}$$

Analogously (5.21) can be transformed to

$$\left(q^3x^2+\frac{s}{q}\right)D^2U_n(x,s,q)+[3]xDU_n(x,s,q)=q^{-n}[n][n+2]U_n(qx,s,q). \qquad (5.23)$$

These equations are of the form considered in [7], (2.2.7). In [7], Section 3.4 for such equations Rodrigues-type formulae have been derived. To transfer these formulae to our special case we need the solution $w(x)$ of the Pearson operator equation $D\left((x^2+s)w(x)\right)=qxw(qx).$

It is easily verified that $w(x)=h\left(-\frac{x^2}{s}\right)$ with

$$h(x)=\sum_{k\geq 0}(-1)^k q^{k^2}\begin{bmatrix}-\frac{1}{2}\\ k\end{bmatrix}_{q^2} x^k=\sum_{k\geq 0}\frac{\left(q;q^2\right)_k}{\left(q^2;q^2\right)_k}x^k=\frac{(qx;q^2)_\infty}{(x;q^2)_\infty}. \qquad (5.24)$$

The last identity follows from the $q-$binomial theorem (2.1).

Obviously $\frac{1}{h(x)}=\sum_{k\geq 0}(-1)^k q^{k^2-k}\begin{bmatrix}\frac{1}{2}\\ k\end{bmatrix}_{q^2} x^k$ is a $q-$analogue of $\sqrt{1-x}.$

Comparing with [7], (3.4.26) we get

**Theorem 5.10 (Rodrigues-type formulae for the q-Chebyshev polynomials)**

$$T_n(x,s,q)=\frac{q^{n(n+1)}}{[1][3]\cdots[2n-1]}\frac{1}{h\left(-\frac{x^2}{s}\right)}D^n\left(h\left(-\frac{x^2}{q^{2n}s}\right)\prod_{k=1}^{n}\left(\frac{x^2}{q^{2k+1}}+\frac{s}{q}\right)\right) \qquad (5.25)$$

*and*

$$U_n(x,s,q)=\frac{[n+1]}{[3][5]\cdots[2n+1]}h\left(\frac{-qx^2}{s}\right)D^n\left(\frac{1}{h\left(\frac{-x^2}{q^{2n-1}s}\right)}\prod_{k=1}^{n}\left(\frac{x^2}{q^{2k-1}}+\frac{s}{q}\right)\right). \qquad (5.26)$$

There are many other interesting formulae. We state some of them and leave the proofs to the reader:

$$T_n(x,s,q) = U_n(x,s,q) - q^n x U_{n-1}(x,s,q), \tag{5.27}$$

$$U_n(x,s,q) = \sum_{k=0}^{n} q^{kn-\binom{k}{2}} x^k \; T_{n-k}(x,s,q), \tag{5.28}$$

$$T_{n+1}(x,s,q) = U_{n+1}(x,s,q) - q^{n+1} x U_n(x,s,q) = x U_n(x,s,q) + q^n s U_{n-1}(x,s,q), \tag{5.29}$$

$$(1+q^n) T_n(x,s,q) = U_n(x,s,q) + q^{2n-1} s U_{n-2}(x,s,q), \tag{5.30}$$

$$U_{2n+1}(x,s,q) = \sum_{k=0}^{n} \left(1+q^{2n+1-2k}\right)(-s)^k q^{4kn+3k-2k^2} T_{2n+1-2k}(x,s,q), \tag{5.31}$$

$$U_{2n}(x,s,q) = \sum_{k=0}^{n-1} \left(1+q^{2n-2k}\right)(-s)^k q^{4kn+k-2k^2} T_{2n-2k}(x,s,q) + (-s)^n q^{2n^2+n}. \tag{5.32}$$

$$T_n(x,q^2 s,q) - T_n(x,s,q) = (q^n - 1) q s U_{n-2}(x,q^2 s,q), \tag{5.33}$$

$$\left(1+q^n\right) T_n(x,s,q) = U_n(x,q^2 s,q) + q s U_{n-2}(x,q^2 s,q), \tag{5.34}$$

$$T_{n+1}(x,s,q) - x T_n(x,s,q) = q^n \left(x^2+s\right) U_{n-1}(x,s,q), \tag{5.35}$$

$$\begin{pmatrix} T_{n+1}(x,s,q) \\ U_n(x,s,q) \end{pmatrix} = \begin{pmatrix} x & q^n\left(x^2+s\right) \\ 1 & q^n x \end{pmatrix} \begin{pmatrix} T_n(x,s,q) \\ U_{n-1}(x,s,q) \end{pmatrix}. \tag{5.36}$$

$$\sum_{n\geq 0} U_n(x,s,q) z^n = \sum_{k\geq 0} q^{\binom{k+1}{2}} z^k x^k \frac{\left(-\frac{sz}{x};q\right)_k}{(xz;q)_{k+1}} \tag{5.37}$$

$$\sum_{n\geq 0} T_n(x,s,q) z^n = \sum_{k\geq 0} q^{\binom{k}{2}} z^k x^k \frac{\left(-\frac{qsz}{x};q\right)_k}{(xz;q)_k}, \tag{5.38}$$

$$U_n(x,s,q)=\det\begin{pmatrix} (1+q)x & qs & 0 & \cdots & 0 & 0 \\ -1 & (1+q^2)x & q^2s & \cdots & 0 & 0 \\ 0 & -1 & (1+q^3)x & \cdots & 0 & 0 \\ \vdots & \vdots & \vdots & \vdots & \vdots & \vdots \\ 0 & 0 & 0 & \cdots & (1+q^{n-1})x & q^{n-1}s \\ 0 & 0 & 0 & \cdots & -1 & (1+q^n)x \end{pmatrix}, \qquad (5.39)$$

$$T_n(x,s,q)=\det\begin{pmatrix} x & qs & 0 & \cdots & 0 & 0 \\ -1 & (1+q)x & q^2s & \cdots & 0 & 0 \\ 0 & -1 & (1+q^2)x & \cdots & 0 & 0 \\ \vdots & \vdots & \vdots & \vdots & \vdots & \vdots \\ 0 & 0 & 0 & \cdots & (1+q^{n-2})x & q^{n-1}s \\ 0 & 0 & 0 & \cdots & -1 & (1+q^{n-1})x \end{pmatrix}. \qquad (5.40)$$

## 6. References

[1 ] W.A. Al-Salam and M.E.H. Ismail, Orthogonal polynomials associated with the Rogers-Ramanujan continued fraction, Pacific J. Math. 104 (1983), 269-283

[2] J. Cigler, q-Fibonacci polynomials, Fib. Quart. 41(2003), 31-40

[3] J. Cigler, Some algebraic aspects of Morse code sequences, DMTCS 6 (2003), 55-68

[4] J. Cigler, q-Fibonacci polynomials and q-Catalan numbers, preprint 2008,

http://homepage.univie.ac.at/johann.cigler/prepr.html

[5] J. Cigler, Some beautiful q-analogues of Fibonacci and Lucas polynomials, arXiv:1104.2699

[6] M.E.H. Ismail, H. Prodinger and D. Stanton, Schur's determinants and partition theorems, Sém. Lothar. Combin. B44a, 2000

[7] R. Koekoek, P.A. Lesky and R.F. Swarttouw, Hypergeometric orthogonal polynomials and their q- analogues, Springer Monographs in Mathematics 2010

[8] M. Schlosser, A non-commutative weight-dependent binomial theorem, arXiv:1106.2112

[9] W. Zudilin, (http://mathoverflow.net/questions/23437/are-the-q-catalan-numbers-q-holonomic)